\documentclass[reqno]{amsart}
\usepackage{graphicx}
\usepackage{amssymb,amsmath}
\usepackage{hyperref}
\usepackage[cmtip,arrow]{xy}
\usepackage{pb-xy}
\usepackage{pb-diagram}
\usepackage[latin1]{inputenc}
\usepackage[T1]{fontenc}
\usepackage{textcomp}
\newtheorem{thm}{Theorem}[section]
\newtheorem{cor}[thm]{Corollary}
\newtheorem{lem}[thm]{Lemma}
\newtheorem{prop}[thm]{Proposition}

\theoremstyle{definition}
\newtheorem{defn}[thm]{Definition}
\theoremstyle{remark}
\newtheorem{rem}[thm]{Remark}
\numberwithin{equation}{section}
\makeindex

\begin{document}
\newcommand{\norm}[1]{\left\Vert#1\right\Vert}
\newcommand{\abs}[1]{\left\vert#1\right\vert}
\newcommand{\set}[1]{\left\{#1\right\}}
\newcommand{\Real}{\mathbb R}
\newcommand{\eps}{\varepsilon}
\newcommand{\To}{\longrightarrow}
\newcommand{\BX}{\mathbf{B}(X)}
\newcommand{\A}{\mathcal{A}}
\newcommand{\gt}{\Gamma_n(T)}
\newcommand{\gn}{\Gamma_n}
\newcommand{\gtt}{\Gamma_n(T)_0}
\newcommand{\ga}{\Gamma_n(A)}
\newcommand{\gm}{\Gamma}
\newcommand{\gaa}{\Gamma _n(A)_0}
\newcommand{\gtn}{\Gamma _n[t_n]}
\newcommand{\gtnn}{\Gamma _n[t_n]_0}
\newcommand{\ggt}{\Gamma_n[t]_0}
\newcommand{\gln}{{\rm Gl}_n(\Z)}
\newcommand{\C}{\mathbb{C}}
\newcommand{\E}{\mathbb{E}}
\newcommand{\F}{\mathbb{F}}
\newcommand{\N}{\mathbb{N}}
\newcommand{\Q}{\mathbb{Q}}
\newcommand{\R}{\mathbb{R}}
\newcommand{\Z}{\mathbb{Z}}

\newcommand{\Hz}{\mathbb{H}}
\newcommand{\Cz}{\C}
\newcommand{\Nz}{\N}
\newcommand{\Qz}{\Q}
\newcommand{\Rz}{\R}
\newcommand{\Zz}{\Z}
\newcommand{\ca}{{\mathcal A}}
\newcommand{\cb}{{\mathcal B}}
\newcommand{\cc}{{\mathcal C}}
\newcommand{\cd}{{\mathcal D}}
\newcommand{\ce}{{\mathcal E}}
\newcommand{\cf}{{\mathcal F}}
\newcommand{\cg}{{\mathcal G}}
\newcommand{\ch}{{\mathcal H}}
\newcommand{\cj}{{\mathcal J}}
\newcommand{\ck}{{\mathcal K}}
\newcommand{\cl}{{\mathcal L}}
\newcommand{\cm}{{\mathcal M}}
\newcommand{\cn}{{\mathcal N}}
\newcommand{\co}{{\mathcal O}}
\newcommand{\cp}{{\mathcal P}}
\newcommand{\cs}{{\mathcal S}}
\newcommand{\cu}{{\mathcal U}}
\newcommand{\cx}{{\mathcal X}}

\newcommand{\sproj}{{\rm Sp_{proj}}(n,\Q)}
\newcommand{\matn}{{\rm M}_n(\Z)}
\newcommand{\drt}{\bigtriangleup_t}
\newcommand{\dr}{\bigtriangleup}
\newcommand{\ii}{{\rm i}}
\newcommand{\comm}{\Upsilon_n}   
\newcommand{\xg}{X_\Gamma}
\newcommand{\fj}[1]{{\mathcal J}_{k,m}(L_1, L_2;#1)}
\newcommand{\fjg}{\fj{\Gamma}}
\newcommand{\oxg}{{\mathcal O}_{\mbox{${\scriptstyle X}_{\scriptscriptstyle \Gamma}$}}}
\newcommand{\calg}{\mathcal G}
\newcommand{\os}{{\mathcal O}_S}
\newcommand{\ox}{{\mathcal O}_X}
\newcommand{\oy}{{\mathcal O}_Y}
\newcommand{\h}{{\rm H}^0}
\newcommand{\xgr}{{X_\Gamma^{{\rm reg}}}}
\newcommand{\ord}{{\rm ord}}

\newcommand{\lto}{\longrightarrow}
\newcommand{\lmto}{\longmapsto}
\newcommand{\mailto}[1]{\href{mailto:#1}{#1}} 
\newcommand{\3}{{\ss}}
\newcommand{\edm}{elementary divisor matrix}
\newcommand{\Ssg}{principal congruence group}
\newcommand{\spec}{{\rm Spec}}
\newcommand{\m}{m}   
\newcommand{\proj}{{\rm Proj}}
\newcommand{\hh}{{\rm H}}
\newcommand{\cech}{\v{C}ech\,}
\newcommand{\pt}{\bullet}      
\newcommand{\htor}{{\rm \underline{Tor}}\,}
\newcommand{\tor}{{\rm {Tor}}}
\newcommand{\goth}[1]{\mathbb{#1}}

\newcommand{\dd}{{\rm d}}
\newcommand{\e}{{\rm e}}
\newcommand{\RG}{\mathbb{R}_{\ge 0}}
\newcommand{\QG}{\mathbb{Q}_{\ge 0}}

\newcommand{\tp}{\tau_\pi}
\newcommand{\vv}[2]{\left .v\right|_{#1}^{#2}\,}
\newcommand{\pr}{{\rm pr}}
\newcommand{\condP}[2]{P(\cdot\vert\,#1)(#2)}   
\newcommand{\rcP}[2]{P^{(#2)}(\cdot\vert\,#1)}  
\newcommand{\w}{\omega}
\newcommand{\pol}[1]{{\cp}[#1,T]}
\newcommand{\stpr}[2]{\{{#1}_t\}_{t\in #2}}

\newcommand{\sigmax}[4]{\cf^{(#1,#2,#3)}_{#4}}
\newcommand{\filt}[3]{\F^{(#1,#2,#3)}}


\newcommand{\optint}[3]{%
            \left[\left.
            \int_{[#1,\tau_#2)}f_s\,\dd
            D_s+1_{\{T<\tau_#2\}}\,h(U_{T-}[#1,#3;#2])
            \right\vert\cf_{[0,#1]}^X
        \right]}

\newcommand{\con}[4]{(#1,#3,#4)}
\newcommand{\f}[2]{\cf^#1_{[0,#2]}}
\newcommand{\sko}{D_E[0,\infty)}
\newcommand{\monoton}{M[0,\infty)}
\newcommand{\bv}{BV[0,\infty)}

\newcommand{\comment}[1]{\marginpar{#1}}

\title[On Skorohod Spaces as Universal Sample Path Spaces]{On Skorohod Spaces as Universal Sample Path Spaces}%
\author{Oliver Delzeith}%
\address{Department of Mathematics, University of Heidelberg\\Im Neuenheimer Feld 280\\
69120 Heidelberg\\Germany}%
\email{\mailto{Delzeith@mathi.uni-heidelberg.de}}%


\commby{Oliver Delzeith}%

\begin{abstract}
The paper presents a factorization theorem for a certain
class of stochastic processes. Skorohod spaces carry the
rich structure of standard Borel spaces and appear to be
suitable universal sample path spaces. We show that, if
$\xi$ is a RCLL stochastic process with values in a
complete separable metric space $E$, any other RCLL
stochastic process $X$ adapted to the filtration induced by
$\xi$ factors through the Skorohod space $D_E[0,\infty)$.
This can be understood as an extension of a stochastic
process to a standard Borel space enjoying nice properties.
Moreover, the trajectories of the factorized stochastic
process defined on $D_E[0,\infty)$ inherit the properties
of being continuous, non-decreasing, and of bounded
variation, resp., from those of $X$.

Considering situations which are invariant under the
factorization procedure, the main theorem is a reduction
tool to assume the underlying measurable space be a
standard Borel space. In an example, we pick the existence
theorem of regular conditional probabilities on standard
Borel spaces to simplify a conditional expectation
appearing in stochastic control problems.
\end{abstract}

\maketitle
\section{Introduction}
\noindent We present a factorization theorem for stochastic
processes applying the theory of standard Borel spaces on
Skorohod spaces.

Let $E$ be a complete separable metric space $E$. The
Skorohod space $\sko$ is defined to be the set of all right
continuous paths $x:[0,\infty)\to E$ having limits on the
left (RCLL). $\sko$ is a metrizable space inheriting the
property to be a complete separable space from $E$
\cite[Thm. 3.5.6, Prop. 3.7.1]{Ethier1986}. Thus, $\sko$
fulfills the defining properties of a standard Borel space
(cf. def. \ref{def:standard Borel space}). The Wiener space
$C_E[0,\infty)$, consisting of continuous paths on
$[0,\infty)$ in E, with the metric of uniform convergence
is contained in $\sko$ as a closed metric subspace, such
that the notion of the Skorohod space generalizes that of
the Wiener space. Just as the Wiener space is the basis of
the stochastic calculus of It{\^o} diffusions and the Brownian
motion (cf., e.g., \cite{Karatzas1997},
\cite{Oksendal2000}, \cite{Ikeda1989}, \cite{Yong1999}),
the present paper shows that the Skorohod space is a
suitable candidate to serve as the fundamental sample path
space for models in financial mathematics and actuarial
sciences. Since stochastic RCLL processes enjoy appropriate
properties to model both payout streams of financial
instruments and the environmental information (e.g., of
capital markets and biometric developments in life
insurance), they  have become of more importance in recent
years. The reader is referred to \cite{Norberg1990},
\cite{Norberg2001}, \cite{Norberg2003},
\cite{Steffensen2001} and \cite{Milbrodt1999a}. The
selected references are not meant to be an exhaustive list.

The paper presents the Skorohod space $\sko$ with its
structure of a complete separable metric space as a
universal sample path space. It will be shown that the
measurable space $(\sko,\cb(\sko))$ carries a natural
filtration $\F^{(E)}=\{\cf_t^{(E )}\}_{t\ge 0}$ together
with a universal stochastic process $\xi^*$ with values in
the complete separable metric space $[0,\infty)\times\sko$
(cf. prop. \ref{prop:universal xi}). If we interpret the
filtration as information getting more precise in the
course of time $t$, then the $\sigma$--algebra $\cf_t^{(E
)}$ carries the information to exactly know the development
of the path $x\in\sko$ up to time $t$. The universal
character of the filtrated measurable space
$$(\sko,\cb(\sko),\F^{(E)})$$ is reflected in the situation
of the factorization theorem \ref{thm:extension of
stochastic processes} which works as follows:

The environmental development of the mathematical model is
driven by a stochastic process
$\xi:[0,\infty)\times\Omega\to E$ on a measurable space
$(\Omega,\cf)$ with RCLL sample paths in a complete
separable metric space $E$. It induces the natural map
\begin{equation}\nonumber
\widetilde{\xi}:\Omega\to\sko\,,\w\mapsto\xi(.,\w)\,.
\end{equation}
The filtration $\{\cf_t^\xi\}_{t\ge 0}$ on $(\Omega,\cf)$
induced by $\xi$ is recovered by $\F^{(E)}$ according to
the relation
\begin{equation}\nonumber
\cf_t^\xi=\widetilde{\xi}^{-1}(\cf_t^{(E )})\,,
\end{equation}
for any $t\ge 0$. Any other observable of the model is
represented by an $\{\cf_t^\xi\}_{t\ge 0}$--adapted
stochastic process $X:[0,\infty)\times\Omega\to E^\prime$
with RCLL sample paths in some separable metric space
$E^\prime$. Note that, because of the assumption on the
sample paths, the notions of adaptedness and progressive
measurability with respect to the filtration
$\{\cf_t^\xi\}$ coincide \cite[Prop. 1.1.13]{Karatzas1997}.
Now, the main result of this paper states that $X$ factors
through the universal sample path space $\sko$ via the
universal map $\widetilde{\xi}$; more precisely, the
factorization theorem provides an $\F^{(E )}$--adapted
stochastic process $X^*:[0,\infty)\times\sko\to E^\prime$
such that
\begin{equation}\label{eq:extension of path I}
X(t,\w)=X^*(t,\widetilde{\xi}(\w))\,,
\end{equation}
for all $(t,\w)\in[0,\infty)\times\Omega$ (cf. the
commutative diagram (\ref{dia:factor intro})).
\begin{equation}\label{dia:factor intro}
\begin{diagram}
\node{[0,\infty)\times\Omega}\arrow{se,t}{X}
\arrow{s,l}{{\rm id}_{[0,\infty)}\times\widetilde{\xi}}\\
\node{[0,\infty)\times\sko}\arrow{e,t,..}{\exists\,
X^*}\node{E^\prime}
\end{diagram}
\end{equation} Rewriting (\ref{eq:extension of path I}) as
\begin{equation}\nonumber
X_t=X^*_t\circ\widetilde{\xi}\,,
\end{equation}
for all $t\ge 0$, shows that the stochastic process
$\{X_t\}_{t\ge 0}$ is extendable to the Skorohod space
$\sko$ via the universal map $\widetilde{\xi}$, basically.
We emphasize that the extension to $\sko$ is not unique, in
general, while we can always extend to the same space
$\sko$, irrespective of the stochastic process $\xi$.

Moreover, it will be shown for special function classes $C$
that, if the sample paths of $X$ belong to $C$, then the
sample paths of the extension $X^*$ can be assumed to be in
class $C$, too. If $C$ denotes the class of continuous
maps, then the result is already known \cite[Thm.
1.2.10]{Yong1999}. We present a proof for the classes $C$
of continuous functions, non-decreasing functions and
functions of bounded variation, resp.\ (cf. cor.
\ref{cor:lifting stochastic process}).

The result of the factorization theorem can easily be
verified if $E$ is a countable discrete space. The present
paper points out that the assertion still holds if the
model states live in, e.g.,  a separable Banach space or a
differentiable manifold.

In an application, we show how the factorization theorem is
the starting point of a proof technique reducing to the
case that the underlying measurable space is a standard
Borel space. After this reduction step is performed, the
rich topological structure of standard Borel spaces implies
nice properties reflected in the measurable structure,
e.g.\ the existence of regular conditional probabilities
with respect to a given sub--$\sigma$--algebra. The given
example appears in the analysis of optimal stochastic
control problems.

\

The paper is organized as follows. In section
\ref{sec:prelims}, we gather preliminary results on
topological spaces and we recall and derive some properties
of standard Borel spaces. In section \ref{sec:filtrated
spaces}, the notion of filtrated Borel spaces is
introduced, and we generalize the Doob-Dynkin lemma to
filtrated Borel spaces. In section \ref{sec:extension}, we
introduce a natural filtration
$\F^{(E)}=\{\cf_t^{(E)}\}_{t\ge 0}$ on the Skorohod space
$\sko$. In the $\sigma$--algebra $\cf_t^{(E)}$, the
information is coded to know the course of the mappings in
$\sko$ up to time $t$. $\cf_t^{(E)}$ emerges as the generic
prototype of the filtration $\{\cf_t^\xi\}_{t\ge 0}$
induced by any stochastic RCLL process $\xi$ with values in
$E$ (prop. \ref{prop:cf^xi}). The proof of the main result
of this paper uses the well known Doob-Dynkin extension
lemma and its generalization to filtrated Borel spaces
developed in section \ref{sec:filtrated spaces}. In section
\ref{sec:application}, we show how to use the factorization
theorem to reduce to the case that the underlying
measurable space is a standard Borel space. The existence
theorem of regular conditional probabilities on standard
Borel spaces yields the simplification of a conditional
expectation which reminds of the strong Markov property of
It{\^o} diffusions (cf., e.g., \cite[Thm.
7.1.2]{Oksendal2000}).

\subsection{Notation} \noindent We fix the following
notational conventions used throughout this article. For a
topological space $E$, the Borel algebra $\cb(E)$ is the
smallest $\sigma$--algebra containing the open sets of $E$.
A measurable map between two measurable spaces is called a
Borel isomorphism iff it is bijective and its inverse is
measurable. If $f:\Omega\to\Omega^\prime$ is a map from a
set $\Omega$ to a measurable space
$(\Omega^\prime,\cf^\prime)$, then $\sigma(f)$ denotes the
smallest $\sigma$--algebra on $\Omega$ such that $f$ is
$\sigma(f)/\cf^\prime$--measurable. For a map
$X:[0,\infty)\times\Omega\to\Omega^\prime$, $\cf_t^X$ is
the $\sigma$--algebra on $\Omega$ generated by
$\sigma(X_s)$, $s\le t$.
If $(\Omega,\cf,\{\cf_t\}_{t\ge 0})$ is a filtrated
measurable space, i.e. $(\Omega,\cf)$ is a measurable space
together with a family of sub--$\sigma$--algebras $\cf_t$
of $\cf$ with $\cf_s\subseteq\cf_t$, $s<t$, then a
stochastic process $X:[0,\infty)\times\Omega\to E$ with
values in a metric space $E$ is called $\{\cf_t\}$--adapted
resp.\ $\{\cf_t\}$--progressively measurable iff, for every
$t\ge 0$, $X_t$ is $\cf_t/\cb(E)$--measurable resp.\ the
induced map $[0,t]\times\Omega\to E$ is
$\cb([0,t])\otimes\cf_t/\cb(E)$--measurable.
For a measurable map
$f:(\Omega,\cf)\to(\Omega^\prime,\cf^\prime)$, the push
forward measure $f_*P$ on $(\Omega^\prime,\cf^\prime)$ of a
measure $P$ on $(\Omega,\cf)$ is defined by $P\circ
f^{-1}$.
A map $x$ from $[0,\infty)$ to a metric space $E$ is called
RCLL iff $x$ is right continuous and has left limits. For a
set $A$, the identity map on $A$ is denoted by ${\rm id}_A$
and the indicator function of $A$ by $1_A$. The notions of
trajectories and sample paths of a stochastic process are
used interchangeably.

\section{Preliminaries}\label{sec:prelims}
\noindent In this section we present results which are
applied in the sequel.

\subsection{Results on topological spaces}

\begin{lem}\label{lem:fixpointset measurable}
Let $E$ be a topological Hausdorff space equipped with the
Borel algebra $\cb(E)$ and $f:E\to E$ a measurable map with
$f\circ f=f$. Then the set of fixpoints of $f$
$${\rm Fix}(f)\stackrel{{\rm def.}}{=}\{e\in E\,\vert\,f(e)=e\}$$
is measurable.
\end{lem}
\proof Since $E$ is Hausdorff, the diagonal
$\Delta_E=\{(e,e)\in E^{\,2}\vert\,e\in E\}$ is a closed,
hence a measurable, subset of $E^{\,2}$. The equation ${\rm
Fix}(f)=(f,{\rm id}_E)^{-1}(\Delta_E)$ together with the
measurable map $(f,{\rm id}_E):E\to E^{\,2}$\index{Identity
${\rm id}_E$} yields the claim. \qed

\begin{lem}\label{lem:fiber product measurable}
Let $(\Omega_i,\cf_i)$, $i=1,2$, be measurable spaces, $E$
a topological Hausdorff space and $f_i:\Omega_i\to E$,
$i=1,2$, measurable maps. Then
\begin{equation}\label{eq:def fiber product}
    \Omega_1\times_E\Omega_2=\{(\w_1,\w_2)\in\Omega_1\times\Omega_2\vert
    \,f_1(\w_1)=f_2(\w_2)\}
\end{equation}
is a measurable subset of $\Omega_1\times\Omega_2$.
\end{lem}
\proof Since $E$ is Hausdorff, the diagonal
$\Delta_E=\{(e,e)\in E^{\,2}\vert\,e\in E\}$ is a closed,
hence Borel measurable, subset of $E^{\,2}$. The map
$f_1\times f_2:\Omega_1\times\Omega_2\to E^{\,2}$ is
measurable. Then the equation
$\Omega_1\times_E\Omega_2=(f_1\times f_2)^{-1}(\Delta_E)$
yields the claim. \qed
\begin{lem}\label{lem:trace Borel algebra}
Let $E$ be a topological space, equipped with the Borel
algebra $\cb(E)$, and $A$ a non-void subset of $E$ equipped
with the relative topology of $E$. Define the
$\sigma$--algebra $\cf^{(A)}=\sigma(\cb(A))$ on $E$. Then:
\begin{enumerate}
\item The Borel algebra $\cb(A)$ of A, equipped with the
relative topology of $E$, coincides with the
trace--$\sigma$--algebra of $\cb(E)$ on $A$:$$\cb(A)=A\cap
\cb(E)\stackrel{{\rm def.}}{=}\{A\cap
B\vert\,B\in\cb(E)\}\,.$$

\item If $(\Omega,\cf)$ is a measurable space and
$f:\Omega\to E$ a map with $f(\Omega)\subset A$, then $f$
is $\cf/\cb(E)$--measurable iff the induced map
$\widetilde{f}:\Omega\to A$ is $\cf/\cb(A)$--measurable iff
$f$ is $\cf/\cf^{(A)}$--measurable.

\item The inclusion $i_A:(A,\cb(A))\hookrightarrow
(E,\cf^{(A)})$ is measurable.
\end{enumerate}
\end{lem}
\begin{rem}
Note that $\cb(A)$ is a $\sigma$--algebra on $E$ iff $A=E$.
\end{rem}
\subsection{Standard Borel spaces} \noindent We recall the
definition of a standard Borel space.
\cite{Parthasarathy1967} represents a general account on
this topic providing the properties of standard Borel
spaces we will apply in this paper.
\begin{defn}\label{def:standard Borel space}
A measurable space $(\Omega,\cf)$ is called \emph{standard
Borel space} iff $(\Omega,\cf)$  is Borel
isomorphic\index{Borel isomorphic} to a complete separable
metric space $(E,\cb(E))$.
\end{defn}
\begin{lem}\label{lem:standard subspaces}
Let $E$ be a standard Borel space and $A$ a measurable
subset of $E$. Then the topological subspace $A$ of $E$ is
a standard Borel space.
\end{lem}
\proof \cite[Thm V.2.2]{Parthasarathy1967}.  \qed
\begin{lem}\label{lem:proj limit/countable index set}
Let $\{E_t\}_{t\ge 0}$ be a family of standard Borel spaces
with a family of surjective measurable maps
$\{\theta_{st}:E_t\to E_s\}_{0\le s<t}$ with
$\theta_{t_1t_3}=\theta_{t_2t_3}\circ\theta_{t_1t_2}$ for
$0\le t_1<t_2<t_3$. Then the projective limit
$\lim\limits_{\leftarrow}E_t$ of $(\{E_t\}_{t\ge
0},\{\theta_{st}\}_{s<t})$ exists in the category of
measurable spaces (together with measurable maps as
morphisms).
\end{lem}
\proof If we restrict the index set to $\N$, the projective
limit $\lim\limits_{\leftarrow}E_n$ exists in the category
of measurable spaces \cite[Thm. V.2.5]{Parthasarathy1967},
which we denote by $E^\prime$. Then it is easy to show that
$E^\prime$ fulfills the universal property of
$\lim\limits_{\leftarrow}E_t$ which completes the proof.
\qed

\section{Filtrated Borel spaces}\label{sec:filtrated spaces}
\subsection{Definition of filtrated Borel spaces}
\begin{defn}\label{def:filtrated space}
A triple $(E,\cf,\{\theta_t\}_{t\ge 0})$ is said to be a
\emph{filtrated Borel space} iff $(E,\cf)$ is a standard
Borel space and $\{\theta_t:E\to E\}_{t\ge 0}$ is a family
of measurable maps satisfying the condition
\begin{equation}\label{eq:formula theta_t}
\theta_t\circ\theta_s=\theta_{s\wedge t}\,,\quad s,t\ge
0\,.
\end{equation}
\end{defn}

An example of filtrated Borel spaces is given by
$(\R,\cb(\R))$ with $\theta_t(s)=s\wedge t$ ($t\ge 0$).
This paper strives for an analysis of Skorohod spaces as
filtrated Borel spaces. We come back to this problem in
section \ref{sec:extension}.

\

Let $(E,\cf,\{\theta_t\}_{t\ge 0})$ be a filtrated Borel
space. The fixpoint set $$E_t={\rm Fix}(\theta_t)=\{e\in
E\vert\,\theta_t(e)=e\}\,,\quad t\ge 0$$ is measurable in
$E$ (lemma \ref{lem:fixpointset measurable}). We equip
$E_t$ with the relative topology of $E$, such that it
becomes a standard Borel space (lemma \ref{lem:standard
subspaces}). The family $\{\theta_t\}_{t\ge 0}$ induces
surjective measurable maps $\widetilde{\theta_t}:E\to E_t$,
$t\ge 0$, with
$\widetilde{\theta_s}={\theta_{st}}\circ\widetilde{\theta_t}$,
$s<t$, where the measurable maps ${\theta_{st}}:E_t\to
E_s$, $s<t$, induced by $\widetilde{\theta_s}$ are
surjective (lemma \ref{lem:trace Borel algebra}).

Using $\{\theta_t\}_{t\ge 0}$, we define a family
$\F^{(E)}=\{\cf_t^{(E)}\}_{t\ge 0}$ of
$\sigma$--sub--algebras of $\cf$ by
\begin{equation}\nonumber
\cf_t^{(E)}=\cf^{(E_t)}=\sigma(\cb(E_t))\,,\quad t\ge 0
\end{equation}
(cf. lemma \ref{lem:trace Borel algebra}.iii). We emphasize
that $\cb(E_t)$ is not a $\sigma$--algebra on $E$, unless
$E_t=E$ which is equivalent to $\theta_t={\rm id}_E$. Since
(\ref{eq:formula theta_t}) implies that $E_s$ is a subset
of $E_t$, for $s<t$, we have $\cb(E_s)=E_s\cap\cb(E_t)$
and, consequently, $\cf_s^{(E)}\subseteq\cf_t^{(E)}$. So
$\F^{(E)}$ is a filtration on $E$. This fact motivates the
notion fixed in definition \ref{def:filtrated space}. The
projective limit $\lim\limits_{\leftarrow}E_t$ of
$(\{E_t\}_{t\ge 0},\{{\theta_{st}}\}_{s<t})$ exists (lemma
\ref{lem:proj limit/countable index set}), and there is a
unique measurable map
${\theta}_E:E\to\lim\limits_{\leftarrow}E_t$ such that the
following diagram is commutative, for any $s<t$:
\begin{equation}\nonumber
\begin{diagram}
\node{E}\arrow{s,l,..}{{\theta}_E}\arrow{se,b}{{\widetilde{\theta}_t}}
\arrow{see,t}{{\widetilde{\theta}_s}}\\
\node{\lim_{\leftarrow}E_t}\arrow{e}\node{E_t}\arrow{e,b}{{\theta_{st}}}\node{E_s}
\end{diagram}
\end{equation}
\begin{defn} A filtrated Borel space $(E,\cf,\{\theta_t\}_{t\ge 0})$ is said to be
\emph{complete} iff the natural measurable map ${\theta}_E$
is bijective.
\end{defn}
If a filtrated Borel space $(E,\cf,\{\theta_t\}_{t\ge 0})$
is complete, ${\theta}_E$ is a Borel isomorphism, by
Kuratowski's theorem \cite[Cor. I.3.3]{Parthasarathy1967}.
Note that, in general, the inclusion $\cup_{t\ge
0}E_t\subseteq E$ is strict, even if
$(E,\cf,\{\theta_t\}_{t\ge 0})$ is supposed to be complete.

\subsection{The Doob-Dynkin lemma on filtrated Borel spaces}
\begin{lem}[Doob-Dynkin: version I]\label{lem:Doob-Dynkin/basic version}
Let $(\Omega,\cf)$, $(\Omega^\prime,\cf^\prime)$ be
measurable spaces and $E$ a standard Borel space. Given any
measurable maps $f:\Omega\to \Omega^\prime$ and
$g:\Omega\to E$, there exists a measurable map
$h:\Omega^\prime\to E$ with $h\circ f=g$ iff the following
condition is satisfied:
\begin{equation}\nonumber
 \sigma(g)\subseteq\sigma(f)\,.
\end{equation}
\end{lem}
\proof \cite[Thm. 1.1.7]{Yong1999}. \qed
\begin{lem}[Doob-Dynkin: version II]\label{lem:Doob-Dynkin/fiber product}
Let $(\Omega,\cf)$ be a measurable space and $E$, $E_i$,
$i=1,2$, standard Borel spaces. Given measurable maps
$f_i:\Omega\to E_i$ and $g_i:E_i\to E$, $i=1,2$, with
$g_1\circ f_1=g_2\circ f_2$, there exists a measurable map
$h:E_1\to E_2$ with $h\circ f_1=f_2$ and $g_2\circ h=g_1$
iff the following condition is satisfied:
\begin{equation}\label{eq:cond. Doob-Dynkin}
 \sigma(f_2)\subseteq\sigma(f_1)\,.\index{$\sigma(f)$}
\end{equation}
\end{lem}
\begin{rem}
We depict the situation of lemma \ref{lem:Doob-Dynkin/fiber
product} in diagram (\ref{dia:map extension}) which can be
completed by an arrow $h$ representing a measurable map
such that the resulting triangles become commutative.
\begin{equation}\label{dia:map extension}
\begin{diagram}
\node[2]{\Omega}\arrow{sw,t}{f_1}\arrow{se,t}{f_2}\\
\node{E_1}\arrow[2]{e,t,..}{\exists
\,h}\arrow{se,b}{g_1}\node[2]{E_2}\arrow{sw,b}{g_2}
\\
\node[2]{E}\\
\end{diagram}
\end{equation}
\end{rem}
\proof The existence of $h$ clearly implies condition
(\ref{eq:cond. Doob-Dynkin}). We show that it is
sufficient, too. The subspace $E_1\times_{E}E_2$ of
$E_1\times E_2$ defined in (\ref{eq:def fiber product}) is
measurable (lemma \ref{lem:fiber product measurable}),
hence a standard Borel space with the relative topology of
$E_1\times E_2$ (lemma \ref{lem:standard subspaces}). The
induced map $(f_1,f_2):\Omega\to E_1\times_{E}E_2$,
$\w\mapsto(f_1(\w),f_2(\w))$ is measurable (lemma
\ref{lem:trace Borel algebra}.ii). We claim that it is
sufficient to construct a measurable map
$\widetilde{h}:E_1\to E_1\times_EE_2$ with
$\widetilde{h}\circ f_1=(f_1,f_2)$. Indeed, in this case
$h={\rm pr}_{E_2}\circ\widetilde{h}$ solves our problem,
where ${\rm pr}_{E_2}:E_1\times_{E}E_2\to E_2$ denotes the
natural measurable projection.
\begin{equation}\label{dia:DoobDynkin/fiber product}
\begin{diagram}
\node{\Omega}\arrow{s,l}{f_1}\arrow{se,b}{(f_1,f_2)}
\arrow{see,t}{f_2}\\
\node{E_1}\arrow{e,b,..}{\widetilde{h}}\node{E_1\times_{E}E_2}
\arrow{e,b}{{\rm pr}_{E_2}}\node{E_2}
\end{diagram}
\end{equation}
So, we can assume to be in the situation depicted in
diagram (\ref{dia:DoobDynkin/fiber product}). Recalling
that $E_1\times_{E}E_2$ is a standard Borel space, version
I of the Doob-Dynkin lemma (lemma
\ref{lem:Doob-Dynkin/basic version}) provides the existence
of a measurable map $\widetilde{h}$ such that the left
triangle in diagram (\ref{dia:DoobDynkin/fiber product})
commutes if $\sigma((f_1,f_2))\subseteq\sigma(f_1)$. The
latter condition is implied by the assumed relation
$\sigma(f_2)\subseteq\sigma(f_1)$. This completes the
proof. \qed
\begin{thm}[Doob-Dynkin: version on filtrated Borel spaces]
\label{thm:Doob-Dynkin/filtrated version} Let
$(\Omega,\cg)$ be a measurable space and
$(E,\cf,\{\theta_t\}_{t\ge 0})$,
$(E^\prime,\cf^\prime,\{\theta_t^\prime\}_{t\ge 0})$
filtrated Borel spaces with the induced filtrations
$\F^{(E)}$ resp.\ $\F^{(E^\prime)}$. Moreover, assume that
\begin{equation}\label{eq:cond measurability}
\cf_t^{(E)}=\widetilde{\theta_t}^{-1}(\cb(E_t))
\quad\text{and}\quad
\widetilde{{\theta_t^\prime}}^{-1}(\cb(E_t^\prime))\subseteq\cf_t^{(E^\prime)}
\,,\end{equation} for any $t\ge 0$, and that
$(E^\prime,\cf^\prime,\{\theta_t^\prime\}_{t\ge 0})$ is
complete.

\noindent Then for any given measurable maps $f:\Omega\to
E$, $f^\prime:\Omega\to E^\prime$, the following assertions
are equivalent:
\begin{enumerate}

\item ${g}^{-1}(\cf_t^{(E^\prime)})\subseteq
f^{-1}(\cf_t^{(E)})$, for each $t\ge 0$;

\item there exists a map $h:E\to E^\prime$ with $h\circ
f=g$ such that
$h^{-1}(\cf_t^{(E^\prime)})\subseteq\cf_t^{(E)}$, for each
$t\ge 0$.
\end{enumerate}
\end{thm}
\begin{rem}
The latter condition of (ii) is equivalent to the
measurability of
$h:(E,\cf_t^{(E)})\to(E^\prime,\cf_t^{(E^\prime)})$, for
each $t\ge 0$.
\end{rem}

\proof Since it is clear that (ii) is sufficient for (i),
we assume that condition (i) holds and perform the proof of
(ii) in three steps. In the first step we fix $t\ge 0$. By
assumption (\ref{eq:cond measurability}), the composed maps
$f_t:(\Omega,\cg)\stackrel{f}{\to}(E,\cf_t^{(E)})\stackrel{\widetilde{\theta_t}}{\to}
(E_t,\cb(E_t))$ and
$g_t:(\Omega,\cg)\stackrel{g}{\to}(E^\prime,\cf_t^{(E^\prime)})
\stackrel{\widetilde{\theta_t^\prime}}{\to}
(E_t^\prime,\cb(E_t^\prime))$ are measurable. Now we want
to construct a measurable map
$h_t:(E_t,\cb(E_t))\to(E_t^\prime,\cb(E_t^\prime))$ with
$h_t\circ f_t=g_t$ (cf. diagram \ref{dia:exist h_t}).
Noting that $E_t^\prime$ is a standard Borel space, we need
only to verify the condition
\begin{equation}\label{eq:extension Doob-Dynkin/proof}
g_t^{-1}(\cb(E_t^\prime))\subseteq f_t^{-1}(\cb(E_t))\,,
\end{equation}
by version I of the Doob-Dynkin lemma (lemma
\ref{lem:Doob-Dynkin/basic version}). But this condition is
clearly implied by the assumptions (\ref{eq:cond
measurability}) and (i).
\begin{equation}\label{dia:exist h_t}
\begin{diagram}
    \node{(\Omega,\cg)}
    \arrow{s,r}{f_t}\arrow{se,t}{g_t}\\
    \node{(E_t,\cb(E_t))}\arrow{e,t,..}{\exists \,h_t}
    \node{(E_t^\prime,\cb(E_t^\prime))}
\end{diagram}
\end{equation}

In the second step, we fix $s,t\ge 0$ with $s<t$. We assume
that there exist measurable maps represented by solid
arrows in the following commutative diagram.
\begin{equation}\nonumber
\begin{diagram}
\node{(\Omega,\cf)}\arrow{se,b}{\widetilde{\theta_t^\prime}\circ
g}\arrow{e,t}{\widetilde{\theta_t}\circ
f}\node{(E_t,\cb(E_t))}\arrow{s,r,..}{h_t}\arrow{e,t}{\theta_{st}}
\node{(E_s,\cb(E_s))}\arrow{s,r}{h_s}\\
\node[2]{(E_t^\prime,\cb(E_t^\prime))}\arrow{e,t}{\theta_{st}^\prime}
\node{(E_s^\prime,\cb(E_s^\prime))}
\end{diagram}
\end{equation}
We want to construct an arrow $h_t$ (represented by the
dotted arrow) such that the resulting triangle and square
become commutative. Again, condition (\ref{eq:extension
Doob-Dynkin/proof}) holds, such that version II of the
Doob-Dynkin lemma \ref{lem:Doob-Dynkin/fiber product}
provides us with a measurable map
$h_t:(E_t,\cb(E_t))\to(E_t^\prime,\cb(E_t^\prime))$ with
the desired properties.

In the third and last step, we build up a measurable map
$h:E\to E^\prime$ with the required properties using the
first two steps and the universal property of the inverse
limit. Since $(E^\prime,\cf^\prime,\{\theta_t^\prime\})$ is
assumed to be complete, $E^\prime$ and
$\lim\limits_{\leftarrow}E_t^\prime$ are Borel isomorphic.
The latter is naturally Borel isomorphic to
$\lim\limits_{\leftarrow}E_n^\prime$ (cf. proof of lemma
\ref{lem:proj limit/countable index set}), such that we can
consider the following diagram (\ref{dia:liftings}) with,
at first, solid arrows only. Going from the right to the
left in diagram (\ref{dia:liftings}), we find completing
dotted arrows representing measurable maps to get a fully
commutative diagram as follows.
\begin{equation}\label{dia:liftings}
\begin{diagram}
\node{\Omega}\arrow{e,t}{f}\arrow{se,b}{g}
\node{E}\arrow{e,t}{\widetilde{\theta}_{n+1}}\arrow{s,l,..}{h}
\node{E_{n+1}}\arrow{e,t}{\theta_{n,n+1}}\arrow{s,l,..}{h_{n+1}}
\node{E_{n}}
\arrow{s,l,..}{h_n}\\
\node[2]{\lim\limits_{\leftarrow}E^\prime_n}
\arrow{e,t}{\widetilde{\theta}_{n+1}^\prime}
\node{E^\prime_{n+1}}\arrow{e,t}{{\theta^\prime}_{n,n+1}}
\node{E^\prime_{n}}
\end{diagram}
\end{equation}
We build a sequence $\{h_n\}_{n\ge 1}$ of measurable maps
$h_n:E_n\to E_n^\prime$ with ${\rm (a)}_n$
$h_{n-1}\circ\theta_{n-1,n}={\theta^\prime}_{n-1,n}\circ
h_n$, for $n\ge 2$, and ${\rm (b)}_n$
$h_n\circ\widetilde{\theta}_n\circ
f=\widetilde{\theta}_n^\prime\circ g$, for $n\ge 1$. To see
this, we apply the first step of the proof to establish
${\rm (b)}_1$ and define $h_n$, $n\ge 1$, inductively with
the second step for ${\rm (a)}_n$ and ${\rm (b)}_n$, $n\ge
2$. Then by the universal property of the inverse limit of
$\{E_n^\prime\}$, the sequence
$\{h_n\circ\widetilde{\theta}_n:E\to E_n^\prime\}_{n\ge 1}$
induces a unique measurable map $h:E\to
\lim\limits_{\leftarrow}E_n^\prime\simeq E^\prime$ with the
desired properties. Thus, the proof is completed. \qed

\section{The factorization theorem of Skorohod
spaces}\label{sec:extension}
\subsection{The Skorohod space $D_E[0,\infty)$} Let $E$ be
a separable complete metric space. We denote the Skorohod
space of $E$ by $\sko$ which is defined as the set of all
maps $x:[0,\infty)\to E$ which are continuous on the right
and have limits on the left (RCLL).\index{RCLL}

$\sko$ becomes a complete separable metrizable space such
that the Borel algebra $\cb(\sko)$ is generated by all
evaluating projections $\pi_t:\sko\to E,x\mapsto x(t)$
($t\ge 0$) \cite[Thm. 3.5.6, Prop. 3.7.1]{Ethier1986}. This
implies a useful criterion of measurability.
\begin{lem}\label{lem:measurablity D_E[0,t]}
If $(\Omega,\cf)$ is a measurable space, a map
$f:\Omega\to\sko$ is measurable iff the composition
$\pi_t\circ f:(\Omega,\cf)\to (E,\cb(E))$ is measurable,
for any $t\ge 0$.
\end{lem}

For $t\ge 0$, we consider the map
\begin{equation}\nonumber
\theta_t:\sko\to\sko, x\mapsto \big(x(.\wedge t):s\mapsto
x(s\wedge t)\big)
\end{equation}
which truncates a map $x$ at $t$ and extends it constantly
on $[t,\infty)$. It holds
$\theta_t\circ\theta_s=\theta_{s\wedge t}$ ($s,t\ge 0$).
$\theta_t$ is measurable, applying lemma
\ref{lem:measurablity D_E[0,t]} and the fact that
$\pi_s\circ\theta_t=\pi_{t\wedge s}$ is measurable, for all
$s,t\ge 0$. Then
$$(\sko,\cb(\sko),\{\theta_t\}_{t\ge 0})$$ is a filtrated Borel space in the sense of section \ref{sec:filtrated spaces}.
Let $D_E[0,t]$ denote the fixpoint set ${\rm
Fix}(\theta_t)$ of $\theta_t$. $D_E[0,t]$ is a measurable
subset of $\sko$ (lemma \ref{lem:fixpointset measurable}).
We have
$$D_E[0,t]=\{x\in\sko\vert x(s)=x(s\wedge t), \forall s\ge
0\}\,.$$ A map $x:[0,\infty)\to E$ is contained in
$D_E[0,t]$ iff $x$ is RCLL and constant on $[t,\infty)$. In
other words, we can naturally identify $D_E[0,t]$ with RCLL
functions $x$ on $[0,t]$ with values in $E$.

As defined in section \ref{sec:filtrated spaces}, the
family $\{\theta_t\}$ induces a filtration
$\F^{(E)}=\{\cf_t^{(E)}\}$ on $\sko$ with
$\cf_t^{(E)}=\sigma(\cb(D_E[0,t]))$.
\begin{lem}\label{lem:generator Borel algebras}
Let $E$ be a complete measurable space. Then it holds, for
any $t\ge 0$,
\begin{equation}\label{eq:gen system/Borel D_E[0,t]}
\cb(D_E[0,t])=\sigma(\{\pi_s\vert_{D_E[0,t]}\vert s\le t\})
\end{equation}
and
\begin{equation}\label{eq:gen system/cf_t}
\cf_t^{(E)}=\widetilde{\theta_t}^{-1}(\cb(D_E[0,t]))=\sigma(\{\pi_s\vert
\,s\le t\})\,.
\end{equation}
In particular, the map
\begin{equation}\label{eq: measurability of theta_t}
\widetilde{\theta_t}:(\sko,\cf_t^{(E)})\to(D_E[0,t],\cb(D_E[0,t]))
\end{equation}
is measurable.
\end{lem}
\proof  Since (\ref{eq:gen system/cf_t}) follows from
(\ref{eq:gen system/Borel D_E[0,t]}), by lemma
\ref{lem:measurablity D_E[0,t]}, and (\ref{eq:
measurability of theta_t}) from (\ref{eq:gen system/cf_t}),
we are reduced to show (\ref{eq:gen system/Borel D_E[0,t]})
which follows from the fact that $\{\pi_t\vert\,t\ge 0\}$
generates $\cb(\sko)$ and $\pi_s^{-1}(B)\cap
D_E[0,t]=\pi_t^{-1}(B)\cap D_E[0,t]$ for $s\ge t$,
$B\subseteq E$. \qed

\subsection{On stochastic processes with RCLL
trajectories}
\noindent The filtrated measurable space
$$(\sko,\cb(\sko),\F^{(E)})$$ features universal properties
which we approach in the remaining parts of this section.

First, we fix some notation. Let $(\Omega,\cf)$ be a
measurable space and $E$ a metric space. Given any
stochastic process $\xi:[0,\infty)\times\Omega\to E$ with
RCLL trajectories in $E$, we define the associated map
\begin{equation}\nonumber
\widetilde{\xi}:\Omega\to \sko\,,\quad
\w\mapsto\xi(.,\w)\,.
\end{equation}
In the sequel, we will use the representations $\xi$,
$\widetilde{\xi}$ and $\{\xi_t\}_{t\ge 0}$ for a stochastic
process interchangeably and analyze their common
properties.
\begin{prop}[Generating property of $\F^{(E)}$]\label{prop:cf^xi}
Let $E$ be a complete separable metric space,
$(\Omega,\cf)$ a measurable space and
$\xi:[0,\infty)\times\Omega\to E$ a stochastic process with
RCLL trajectories. Then
\begin{equation}\label{eq:cf^xi}
\cf_t^\xi=\widetilde{\xi}^{-1}(\cf_t^{(E)})\quad(t\ge 0)\,.
\end{equation}
\end{prop}
\proof Due to $\xi_s=\pi_s\circ\widetilde{\xi}$, it holds
$\xi^{-1}_s(B)=\widetilde{\xi}^{-1}(\pi_s^{-1}(B))$ for any
subset $B\subset E$, $s\ge 0$. By (\ref{eq:gen system/Borel
D_E[0,t]}), the generating system
$\{\xi^{-1}_s(B)\vert\,s\le t,\,B\in\cb(E)\}$ of
$\cf_t^\xi$ generates $\widetilde{\xi}^{-1}(\cf_t^{(E)})$,
too, which implies the claim. \qed

\begin{lem}\label{lem:measuribility to sko}
Let $T>0$, $(\Omega,\cf)$ be a measurable space and $E$ be
a complete separable metric space. For any stochastic
process $\xi:[0,T]\times \Omega\to E$ with RCLL
trajectories, the following assertions are equivalent:
\begin{enumerate}
    \item $\widetilde{\xi}$ is $\cf/\cb(D_E[0,T])$--measurable;

    \item $\xi_t$ is $\cf/\cb(E)$--measurable, for each $t\in[0,T]$;

    \item ${\xi}$ is $\cb[0,T]\otimes\cf/\cb(E)$--measurable.

\end{enumerate}
\end{lem}
\proof The two facts that we have
$\xi_t=\pi_t\circ\widetilde{\xi}$ and that the Borel
algebra of $D_E[0,T]$ is generated by
$\{\pi_t\vert_{D_E[0,T]}\vert\,t\in[0,T]\}$ (lemma
\ref{lem:generator Borel algebras}) imply the equivalence
of (i) and (ii). The assumption that $\xi$ has right
continuous trajectories yields the implication from (ii) to
(iii), by \cite[Prop. 1.13]{Karatzas1997}. Note that the
reference provides the proof for this standard result if
$E=\R^n$. The argument given therein can be immediately
generalized to any complete metric space $E$ as it is. The
reverse implication is trivial. So, we have shown the
equivalence of all stated conditions. \qed

\begin{lem}\label{lem:progressively measurable}
Let $(\Omega,\cf,\{\cf_t\}_{t\ge 0})$ be a filtrated
measurable space, $E$ a complete separable metric space and
$\xi:[0,\infty)\times\Omega\to E$ a stochastic process with
RCLL trajectories. Then the following assertions are
equivalent:
\begin{enumerate}
    \item  $\xi$ is $\{\cf_t\}_{t\ge 0}$--progressively
    measurable;

    \item for any $t\ge 0$, $[0,t]\times\Omega\to E$, $(s,\w)\mapsto \xi(s,\w)$
    is $\cb([0,t])\otimes\cf_{t}/\cb(E)$--measurable;

    \item for any $t\ge 0$, $\Omega\to D_E[0,t]$, $\w\mapsto
    \xi(.,\w)\vert[0,t]$ is
    $\cf_t/\cb(D_E[0,t])$--measurable;

    \item for any $t\ge 0$, $\Omega\to D_E[0,\infty)$,
    $\w\mapsto\xi(.,\w)\vert[0,t]$ is
    $\cf_t/\cf_t^{(E)}$--measurable.
\end{enumerate}
\end{lem}
\begin{rem}
Note that the set-theoretic maps in (iii) and (iv) coincide
with $\widetilde{\theta_t}\circ\widetilde{\xi}$ and
$\theta_t\circ \widetilde{\xi}$, respectively.
\end{rem}

\proof The stated equivalence follows from lemma
\ref{lem:measuribility to sko} and the fact that the
measurable map $\widetilde{\xi}:\Omega\to\sko$ having
values in $D_E[0,t]$ factorizes through the
$\cf_t/\cf_t^{(E)}$--measurable inclusion
$i_t:D_E[0,t]\hookrightarrow\sko$ (cf. lemma \ref{lem:trace
Borel algebra}). \qed

\begin{prop}\label{prop:universal xi}
Let $(\Omega,\cf)$ be a measurable space and $E$ a complete
separable metric space. Then
\begin{equation}\nonumber
\Omega_E^*\stackrel{{\rm def.}}{=}\sko\quad\text{and}\quad
E^*\stackrel{{\rm def.}}{=}[0,\infty)\times\sko\,.
\end{equation}
are complete separable metric spaces and the stochastic
process
\begin{equation}\nonumber
\xi^*:[0,\infty)\times \Omega_E^*\to
E^*\,,(t,x)\mapsto(t,\theta_t(x))
\end{equation}
has RCLL trajectories. It holds, for any $t\ge 0$,
\begin{equation}\label{eq:xi^* submersive}
\widetilde{\xi^*}^{-1}(\cf_t^{(E^*)})=\cf_t^{(E)}\,.
\end{equation}
\end{prop}
\begin{rem}
Note that $[0,\infty)$ is equipped with the family
$\theta_t(s)=s\wedge t$ ($s,t\ge 0$) to establish a
filtrated Borel space structure on the product space $E^*$.
\end{rem}
\proof Since $\sko$ is a complete separable metric space,
so are $\Omega_E^*$ and $E^*$. Let denote $d^*$ a complete
metric on $\Omega_E^*=\sko$ (cf. \cite[Chap.
3.5]{Ethier1986} for details on $d^*$). Now, we show that
$\xi^*$ has RCLL trajectories in $E^*$. For this purpose,
it suffices to prove that $t\mapsto\theta_t(x)$ is
right-continuous and has left limits, for any fixed
$x\in\sko$. By right-continuity of $x$, $\theta_{t_n}x$
converges uniformly to $\theta_t(x)$, as $t_n\downarrow t$,
which implies that $\theta_{t_n}(x)$ converges to
$\theta_t(x)$ under the metric $d^*$. Similarly, it can be
shown that $\theta_{t_n}x$ converges uniformly to
$x^{(t)}=x\,1_{[0,t)}+x(t-)\,1_{[t,\infty)}$ with
$x(t-)=\lim\limits_{s\uparrow t}x(s)$ in $\sko$, as
$t_n\uparrow t$, hence $\theta_{t_n}x\to\theta_tx$ under
$d^*$. This implies that $t\mapsto\theta_tx$ is RCLL. Next,
we prove relation (\ref{eq:xi^* submersive}). By prop.
\ref{prop:cf^xi}, we have
$\widetilde{\xi^*}^{-1}(\cf_t^{(E^*)})=\cf_t^{\xi^*}$, such
that we are reduced to show
\begin{equation}\nonumber
\cf_t^{\xi^*}=\cf_t^{(E)}\,,
\end{equation}
for all $t\ge 0$. To see this, we present a common
generating system of the two given $\sigma$--algebras by
\begin{equation}\nonumber
\{\xi_s^*\vert s\le t\}=\{\theta_s\vert s\le
t\}\stackrel{{\rm lem.} \ref{lem:measurablity D_E[0,t]}}{=}
\{\pi_{s^\prime}\circ\theta_s\vert s^\prime\le s\le t\}
=\{\pi_s\vert s\le t\}\,,
\end{equation}
where we used $\xi_s^*(x)=(s,\theta_s(x))$ and
$\pi_{s^\prime}\circ\theta_s=\pi_{s^\prime}$, for
$s^\prime\le s$. This completes the proof. \qed

\subsection{The factorization theorem}

After the preparation work, we are ready to prove the main
result of this paper.
\begin{thm}\label{thm:extension of stochastic processes}
Let $(\Omega,\cf)$ be a measurable space, $E$ a complete
separable metric space. Define the standard Borel space
$\Omega^*_E=\sko$. Then the filtrated measurable space
$$(\Omega^*_E,\cb(\Omega^*_E),\{\cf_t^{(E)}\}_{t\ge 0})$$
is universal in the following sense:

\noindent Let $\xi:[0,\infty)\times\Omega\to E$ be a
stochastic process with RCLL trajectories inducing the
filtration $\{\cf_t^\xi\}_{t\ge 0}$ on $(\Omega,\cf)$ and
the map $\widetilde{\xi}:\Omega\to\Omega^*_E$. Then for any
separable metric space $E^\prime$ and any
$\{\cf_t^\xi\}_{t\ge 0}$--progressively measurable
stochastic process $X:[0,\infty)\times\Omega\to E^\prime$
with RCLL trajectories, there exists an
$\{\cf_t^{(E)}\}_{t\ge 0}$--progressively measurable
stochastic process $X^*:[0,\infty)\times\Omega^*_E\to
E^\prime$ with RCLL trajectories such that the following
diagram
\begin{equation}\label{dia:extension process}
\begin{diagram}
\node{[0,\infty)\times\Omega} \arrow{see,t}{{X}}
\arrow{s,l}{{\rm id}_{[0,\infty)}\times\widetilde{\xi}}\\
\node{[0,\infty)\times\Omega_E^*}\arrow[2]{e,t,..}{\exists
\,{X^*}}\node[2]{E^\prime}
\end{diagram}
\end{equation}
is commutative, i.e. it holds, for
$(t,\w)\in[0,\infty)\times\Omega$,
$$X(t,\w)=X^*(t,\widetilde{\xi}(\w))=X^*(t,\xi(.,\w))\,.$$
\end{thm}
\begin{rem}
Note that thm. \ref{thm:extension of stochastic processes}
does not claim that $X^*$ in diagram (\ref{dia:extension
process}) is uniquely determined. But the theorem
guarantees that the Skorohod space $\Omega^*_E=\sko$
provides enough information about any stochastic process
with values in $E$ such that any process $X$ with values in
any other separable metric space $E^\prime$ can be extended
from $\Omega$ to $\Omega^*_E$.
\end{rem}

\proof To start with, we note that we can assume $E^\prime$
to be complete, without loss of generality, by taking the
completion of $E^\prime$ if necessary. We will apply the
version of the Doob-Dynkin lemma for filtrated Borel spaces
(prop. \ref{thm:Doob-Dynkin/filtrated version}). By prop.
\ref{prop:universal xi}, $\xi^*$ has RCLL trajectories and
induces the well-defined map
$\widetilde{\xi^*}:\Omega_E^*\to D_{E^*}[0,\infty)$. Define
$\xi^\#=\widetilde{\xi^*}\circ\widetilde{\xi}:\Omega\to\Omega_E^*\to
D_{E^*}[0,\infty)$. We consider the following diagram
(\ref{dia:extension I}) with, at first, solid arrows only.
\begin{equation}\label{dia:extension I}
\begin{diagram}
\node{\Omega}\arrow{se,t}{\widetilde{X}}\arrow{s,l}{{\xi}^\#}\\
\node{D_{E^*}[0,\infty)}\arrow{e,t,..}{Y}\node{D_{E^\prime}[0,\infty)}
\end{diagram}
\end{equation}
By assumption and lemma \ref{lem:progressively measurable},
it holds
\begin{equation}\nonumber
\widetilde{X}^{-1}(\cf_t^{(E^\prime)})\subseteq\cf_t^{\xi}
\end{equation}
and
\begin{equation}\nonumber
{\xi^\#}^{-1}(\cf_t^{(E^*)})=
\widetilde{\xi}^{-1}(\widetilde{\xi^*}^{-1}(\cf_t^{(E^*)}))
\stackrel{(\ref{eq:xi^* submersive})}{=}
\widetilde{\xi}^{-1}(\cf_t^{(E)})
\stackrel{(\ref{eq:cf^xi})}{=} \cf_t^\xi\,,
\end{equation}
for all $t\ge 0$, implying the Doob-Dynkin condition
\begin{equation}\nonumber
\widetilde{X}^{-1}(\cf_t^{(E^\prime)})\subseteq
{{{\xi^\#}}}^{-1}(\cf_t^{(E^*)})\,,
\end{equation} for any $t\ge 0$.
Since the additional condition (\ref{eq:cond
measurability}) of thm. \ref{thm:Doob-Dynkin/filtrated
version} is given by lemma \ref{lem:generator Borel
algebras}, we have verified all assumptions of thm.
\ref{thm:Doob-Dynkin/filtrated version} such that there
exists a map $Y:D_{E^*}[0,\infty)\to
D_{E^\prime}[0,\infty)$ with $Y\circ\xi^\#=X$ and
\begin{equation}\label{eq:DD cond}
Y^{-1}(\cf_t^{(E^\prime)})\subseteq\cf_t^{(E^*)}\,,
\end{equation}
for all $t\ge 0$. Then
\begin{equation}\nonumber
X^*:[0,\infty)\times\Omega_E^*\to E^\prime\,,(t,x)\mapsto
(Y\circ\widetilde{\xi})(x)(t)\,
\end{equation}
defines a stochastic process with
$X(t,\w)=X^*(t,\widetilde{\xi}(\w))$ which is
${\cf_t^{(E)}}$--progressively measurable, by lemma
\ref{lem:progressively measurable}, (\ref{eq:xi^*
submersive}) and (\ref{eq:DD cond}). This completes the
proof of the theorem. \qed

We want to clarify the importance of thm.
\ref{thm:extension of stochastic processes} in terms of
universal mappings.
\begin{defn}
Let $(\Omega,\cf,\{\cf_t\})$ be a filtrated measurable
space and $E^\prime$ a metric space. Then
$$D(\Omega,\{\cf_t\};E^\prime)$$ denotes the set of all
$\{\cf_t\}$--progressively measurable stochastic processes
$X:[0,\infty)\times\Omega\to E^\prime$ with RCLL
trajectories.
\end{defn}
\begin{rem}\label{rem:adapted/progressively measurable}
Note that a stochastic process $X[0,\infty)\times\Omega\to
E^\prime$ with RCLL trajectories is in
$D(\Omega,\{\cf_t\};E^\prime)$ iff $X$ is
$\{\cf_t\}$-adapted, by \cite[Prop. 1.1.13]{Karatzas1997}.
\end{rem}
\begin{cor}
In the situation of thm. \ref{thm:extension of stochastic
processes}, the map
\begin{equation}\nonumber
\left\{
\begin{array}{ccl}
D(\Omega^*_E,\{\cf_t^{(E)}\};E^\prime)&\longrightarrow&D(\Omega,\{\cf_t^\xi\};E^\prime)\\
X^*&\longmapsto&X^*\circ({\rm
id}_{[0,\infty)}\times\widetilde{\xi})
\end{array}\right.
\end{equation}
induced by the universal stochastic process
$${\rm
id}_{[0,\infty)}\times\widetilde{\xi}:[0,\infty)\times\Omega\to[0,\infty)\times\sko$$
with RCLL trajectories associated with $\xi$ is a
surjection.
\end{cor}
\begin{prop}\label{prop:monoton/bounded var=standard
spaces} Let $\monoton$ resp.\ $\bv$ denote the set of
functions in $D_\R[0,\infty)$ which are non-decreasing
resp.\  of bounded variation. Then $\monoton$ and $\bv$ are
measurable in $D_\R[0,\infty)$. Moreover, under the
filtration structure $\{\theta_t^{(E)}\}_{t\ge 0}$
inherited from $D_\R[0,\infty)$, $\monoton$ and $\bv$
become complete filtrated Borel spaces.
\end{prop}
\proof Firstly, we proof that $\monoton$ is measurable in
$D_\R[0,\infty)$. Define $I=\{(s,t)\in\RG\vert\,s<t\}$. The
set
$$A(s,t)=\{x\in\sko\vert\,x(s)\le x(t)\}\quad ((s,t)\in I)$$
is measurable. Indeed, note that the set
$B=\{(a,b)\in\R^2\vert\,a\le b\}$ is a measurable subset of
$\R^2$, the product map
$\pi_s\times\pi_t:D_\R[0,\infty)^2\to\R^2$ resp.\ the
diagonal map $\Delta:D_\R[0,\infty)\to D_\R[0,\infty)^2$,
$x\mapsto(x,x)$ are measurable and it holds
$$A(s,t)=((\pi_s\times\pi_t)\circ\Delta)^{-1}(B)\,.$$
We have
\begin{equation}\nonumber
    \monoton=\bigcap_{(s,t)\in I}A(s,t)=
    \bigcap_{(p,q)\in I\cap\Q^2}A(p,q)\,,
\end{equation}
where the second equality is implied by the right
continuity of any function $x$ in $D_\R[0,\infty)$, such
that we can represent $\monoton$ as a countable
intersection of measurable sets, hence $\monoton$ is
measurable in $D_\R[0,\infty)$.

Secondly, we show the measurability assertion on $\bv$.
Define $\monoton_0$ to be the set of functions
$x\in\monoton$ with $x(0)=0$. $\monoton_0$ is measurable,
since $\monoton_0=\monoton\cap\pi_0^{-1}(\{0\})$, hence a
standard Borel space, by the first step and lemma
\ref{lem:standard subspaces}. Now, we apply the well known
fact that any function $x\in\bv$ can be uniquely decomposed
in $x=a+x_1-x_2$ with $x_1,x_2\in\monoton_0$, $a\in\R$,
i.e.
$$\Phi:\monoton_0^2\times\R\to\bv\,,(x_1,x_2,a)\mapsto a+x_1-x_2$$
is a bijection (see, e.g., \cite[Satz 14.4]{Meintrup2004}).
$\Phi$ is a measurable map considered as a map to
$D_\R[0,\infty)$, by lemma \ref{lem:measurablity D_E[0,t]}.
By \cite[Thm. I.3.9]{Parthasarathy1967}, the image of
$\Phi$ is measurable in $D_\R[0,\infty)$.

As $\theta_t^{(E)}$ leaves both $\monoton$ and $\bv$
invariant, the additional claim on the completeness as
standard Borel spaces is easily verified, since, for any
function $x$, the defining properties of $\monoton$ and
$\bv$ can be checked on the compact intervals $[0,t]$,
$t\ge 0$. \qed

\begin{lem}\label{lem:C closed in D} Let $E$ be a complete
separable metric space. Then the space $C_E[0,\infty)$ of
continuous maps $x:[0,\infty)\to E$ with the metric of
uniform convergence is a closed subspace of $\sko$ with the
relative topology. Moreover, together with the filtration
structure $\{\theta_t^{(E)}\}_{t\ge 0}$ inherited from
$\sko$, $C_E[0,\infty)$ becomes a complete filtrated Borel
space.
\end{lem}
\proof \cite[Lemma VII.6.8]{Parthasarathy1967}. \qed

\begin{cor}\label{cor:lifting stochastic process}
Suppose that the situation of thm. \ref{thm:extension of
stochastic processes} holds and let $E$ be a complete
separable metric space. Then, if the trajectories of $X$
belong to $\monoton$, $\bv$ and $C_E[0,\infty)$, resp., the
trajectories of $X^*$ can be assumed to have the same
property.
\end{cor}
\noindent The \emph{proof} follows immediately with thm.
\ref{thm:extension of stochastic processes}, prop.
\ref{prop:monoton/bounded var=standard spaces} and lemma
\ref{lem:C closed in D}. \qed

\section{An application of the factorization theorem}\label{sec:application} \noindent We want to show
an application of thm. \ref{thm:extension of stochastic
processes} resp.\ cor. \ref{cor:lifting stochastic process}
providing the comfortable situation to assume that the
underlying probability space is a standard Borel space.
This justifies the technical work we have gone through in
the previous sections.
\begin{lem}
Let $\xi:\Omega\to\Omega^\prime$ be a measurable map from
the measure space $(\Omega,\cf,P)$ to the measurable space
$(\Omega^\prime,\cf^\prime)$ and $\cg^\prime$ a
sub--$\sigma$--algebra of $\cf^\prime$. Define
$\cg=\xi^{-1}(\cg^\prime)$. \newline If
$f:\Omega^\prime\to\R$ is an
$\cf^\prime/\cb(\R)$--measurable map, then it holds
\begin{equation}\label{eq:f_* cond expectation}
    \E^{\xi_*P}[f\vert\cg^\prime]\circ\xi=\E^P[f\circ\xi\vert\cg]\quad
    \text{$P\vert_\cg$--a.s.}
\end{equation}
\end{lem}
\proof Checking that
$\E^{\xi_*P}[f\vert\cg^\prime]\circ\xi$ satisfies the
defining Radon-Nikod{\'y}m equation of the conditional
expectation of $\E^P[f\circ\xi\vert\cg]$ gives the proof.
\qed

\

As an example to show how the factorization theorem
\ref{thm:extension of stochastic processes} works as a
reduction tool, we consider a situation appearing in
stochastic control problems.
\begin{thm}
Let $(\Omega,\cf,P)$ be a measure space,
$\xi:[0,\infty)\times\Omega\to E$ a stochastic process with
RCLL trajectories in a complete separable metric space $E$
and $X:[0,\infty)\times\Omega\to\R$ an
$\{\cf_t^\xi\}$--progressively measurable process with RCLL
trajectories of bounded variation. For the initial value
problem
\begin{equation}\nonumber
(S)\quad
\left\{
\begin{array}{lll}
\dd Y_t=\e^{-\alpha t}\dd X_t&;&t\in[t_0,\infty)\\
Y_{t_0}=Z
\end{array} \right.
\end{equation}
with an $\cf_{t_0}^\xi$--measurable random variable
$Z:\Omega\to \R$ and a constant $\alpha>0$, denote the
(pathwise uniquely determined) solution by $Y_t^{(Z)}$.
Then for $t\ge t_0$ and any bounded measurable map
$f:\R\to\R$, it holds
\begin{equation}\label{eq:quasi Markov property}
\E^P[f\circ Y_t^{(Z)}\vert\,\cf_{t_0}^\xi]= \E^P[f\circ
Y_t^{(z)}\vert\,\cf_{t_0}^\xi]\big\vert_{z=Z}\,.
\end{equation}
\end{thm}
\begin{rem}
The solution of (S) can be represented as a pathwise
Lebesgue-Stieltjes integral. Then the stochastic process
$Y_t^{(Z)}:[0,\infty)\times\Omega\to\R$ is
$\{\cf_t^\xi\}$--progressively measurable; this is because
it has right continuous trajectories and is
$\{\cf_t^\xi\}$-adapted, such that the claim follows with
remark \ref{rem:adapted/progressively measurable}. In
particular, the conditional expectations in (\ref{eq:quasi
Markov property}) are well-defined.
\end{rem}
\proof In the first step of the proof, we show that we can
confine ourselves to the case that $\Omega$ can be assumed
to be a standard Borel space. To see this, we apply the
basic version of the Doob-Dynkin lemma (lemma
\ref{lem:Doob-Dynkin/basic version}) and cor.
\ref{cor:lifting stochastic process} to extend $Z$ and $X$
to the standard Borel space $\Omega_E^*$, i.e. there are an
$\cf_{t_0}^{(E)}$--measurable random variable
$Z^*:\Omega_E^*\to \R$ and an
$\{\cf_t^{(E)}\}$--progressively measurable process
$X^*:[0,\infty)\times\Omega_E^*\to \R$ with RCLL
trajectories of bounded variation such that
$$Z=Z^*\circ\widetilde{\xi}
\quad\text{and}\quad X_t=X^*_t\circ \widetilde{\xi}\,,$$
for all $t\ge 0$. The induced initial value problem
\begin{equation}\nonumber(S^*)\quad
\left\{
\begin{array}{lll}
\dd Y_t^*=\e^{-\alpha t}\dd X_t^*&;&t\in[t_0,\infty)\\
Y_{t_0}^*=Z^*
\end{array} \right.
\end{equation}
has a (pathwise unique) solution which is denoted by
$Y_t^{*(Z^*)}$. The pathwise uniqueness of the solutions of
(S) implies, for any $t\ge 0$,
$$Y_t^{(Z)}=Y_t^{*(Z^*)}\circ \widetilde{\xi}\,.$$ Thus, we
get
\begin{eqnarray*}
\E^P[f(Y_t^{(Z)}\vert\cf_{t_0}^\xi] &=& \E^P[f\circ
Y_t^{*(Z^*)}\circ\widetilde{\xi}\vert
\widetilde{\xi}^{-1}(\cf_{t_0}^{(E)})]\\
&\stackrel{(\ref{eq:f_* cond
expectation})}{=}&\E^{\widetilde{\xi}_*P}[f\circ
Y_t^{*(Z^*)}\vert\cf_{t_0}^{(E)}]\circ\widetilde{\xi}
\end{eqnarray*}
and \begin{eqnarray*} \E^P[f(Y_t^{(z)}\vert\cf_{t_0}^\xi]
&=& \E^P[f\circ Y_t^{*(z)}\circ\widetilde{\xi}\vert
\widetilde{\xi}^{-1}(\cf_{t_0}^{(E)})]\\
&\stackrel{(\ref{eq:f_* cond
expectation})}{=}&\E^{\widetilde{\xi}_*P}[f\circ
Y_t^{*(z)}\vert\cf_{t_0}^{(E)}]\circ\widetilde{\xi}\,,
\end{eqnarray*}
such that, for proving (\ref{eq:quasi Markov property}), it
suffices to show
$$\E^{\widetilde{\xi}_*P}[f\circ
Y_t^{*(Z^*)}\vert\cf_{t_0}^{(E)}]=\E^{\widetilde{\xi}_*P}[f\circ
Y_t^{*(z)}\vert\cf_{t_0}^{(E)}]\big\vert_{z=Z^*}\,.$$

To begin with the second step, we note that the result of
the first step allows us to assume $\Omega$ to be a
standard Borel space, without loss of generality. Now, we
apply the fact that, on standard Borel spaces, regular
conditional probabilities with respect to any
sub--$\sigma$--algebra exist (cf. \cite[Chap.
5.3.C]{Karatzas1997} and the references therein, or
\cite[Thm. 3.1]{Leao2004}, in terms of Radon spaces). Set
$\cg=\cf_{t_0}^\xi$. Denote by $P(\cdot\vert\cg)$ the
regular conditional probability on $(\Omega,\cf)$ given
$\cg$. We prove (\ref{eq:quasi Markov property}) by
evaluating the equation pointwise at an arbitrary argument
$\w_0\in\Omega$. By \cite[Thm. 5.3.18]{Karatzas1997}, $Z$
is constant almost surely under the measure
$P(\cdot\vert\cg)(\w_0)$; more precisely we have
$$Z(\w)=Z(\w_0)\quad\text{$P(\cdot\vert\cg)(\w_0)$--almost all $\w\in\Omega$.}$$
This yields, for any $\w_0\in\Omega$,
\begin{eqnarray*}
E^P[f\circ Y_t^{(Z)}\vert\cg](\w_0)&=&
E^{P(\cdot\vert\cg)(\w_0)}[f\circ Y_t^{(Z)}]\\
&=&E^{P(\cdot\vert\cg)(\w_0)}[f\circ Y_t^{(Z(\w_0))}]\\
&=&E^{P(\cdot\vert\cg)(\w_0)}[f\circ Y_t^{(z)}]\big\vert_{z=Z(\w_0)}\\
&=&E^P[f\circ
Y_t^{(z)}\vert\cg]\big\vert_{z=Z(\w_0)}(\w_0)\,.
\end{eqnarray*}
This completes the proof of (\ref{eq:quasi Markov
property}). \qed

\providecommand{\bysame}{\leavevmode\hbox to3em{\hrulefill}\thinspace}
\providecommand{\MR}{\relax\ifhmode\unskip\space\fi MR }
\providecommand{\MRhref}[2]{%
  \href{http://www.ams.org/mathscinet-getitem?mr=#1}{#2}
}
\providecommand{\href}[2]{#2}

\end{document}